\documentclass[titlepage,11pt]{article}
\usepackage{latexsym}
\usepackage{amsfonts}
\usepackage{amsmath,amsthm}
\newcommand\blackslug{\hbox{\hskip 1pt \vrule width 4pt height 8pt depth 1.5pt
        \hskip 1pt}}
\newcommand\bbox{\hfill \quad \blackslug \bigbreak}
\def\DD{\hbox{-}}

\def\LL{,\ldots,}

\usepackage[utf8]{inputenc}
\usepackage[T1]{fontenc}

\usepackage{url} 
\usepackage[hypertexnames=false]{hyperref}
\usepackage[noabbrev,capitalise]{cleveref}

\usepackage{enumitem}

\newtheorem{thm}{}[section]

\newcommand{\Proof}{\noindent{\bf Proof.}\ \ }

\newtheoremstyle{claimstyle} 
    {}                    
    {}                    
    {\itshape}                   
    {}                           
    {(}                   
    {)}                          
    {5pt}                       
    {}  
\theoremstyle{claimstyle}
\newtheorem{claim}{}
\creflabelformat{claim}{#2(\textup{#1})#3}

\usepackage{chngcntr}
\counterwithin*{claim}{section}

\title{Cops and robbers on $P_5$-free graphs}
\author{Maria Chudnovsky\thanks{Supported by NSF DMS-EPSRC grant DMS-2120644 and AFOSR grant FA9550-22-1-0083.}\\
Princeton University, Princeton, NJ 08544
\\
\\
Sergey Norin\thanks{Supported by an NSERC Discovery Grant. Support\'e par le Programme de subventions \`a la recherche du CRSNG.}
\\Department of Mathematics and Statistics, McGill University, Montr\'eal, QC, Canada
\\
\\
Paul Seymour\thanks{Supported by AFOSR grant
FA9550-22-1-0234, and by NSF grant DMS-2154169.}\\
Princeton University, Princeton, NJ 08544
\\
\\
J\'er\'emie Turcotte\thanks{Supported by an NSERC CGS D scholarship. Supporté par une bourse BESC D du CRSNG.}\\
Department of Mathematics and Statistics, McGill University, Montr\'eal, QC, Canada}

\date{December 19, 2022; revised January 25, 2023}

\begin{document}
\maketitle
\begin{abstract}
We prove that every connected $P_5$-free graph has cop number at most two, solving a conjecture of Sivaraman. In order to do so, we first prove that every connected $P_5$-free graph $G$ with independence number at least three contains a three-vertex induced path with vertices $a \hbox{-} b \hbox{-} c$ in order, such that every neighbour of $c$ is also adjacent to one of $a,b$.
\end{abstract}

\section{Introduction}
We denote the $t$-vertex path by $P_t$. There are a number of well-known open questions about $P_5$-free graphs (a graph $G$ is {\em $H$-free} if no induced subgraph of $G$ is isomorphic to $H$, and $|G|$ denotes the number of vertices of $G$). For instance:
\begin{itemize} 
	\item the Erd\H{o}s-Hajnal \cite{erdos_ramsey-type_1989} conjecture implies that for some $c>0$, every $P_5$-free graph $G$ has a clique or stable set of size at least $|G|^c$;
	\item a conjecture of Esperet~\cite{esperet_graph_2017} implies that for some $c>0$, every $P_5$-free graph $G$ has chromatic number at most $\omega(G)^c$, where $\omega(G)$ is the clique number of $G$;
	\item a conjecture of Ho\`ang et al.~\cite{hoang_deciding_2010} says that for some $c>0$, there is a function $f$ such that for every $k$ there is an algorithm deciding whether a $P_5$-free graph $G$ is $k$-colourable in time $f(k)|G|^c$.
\end{itemize}

In this paper, we study another conjecture on $P_5$-free graphs, which concerns the game of cops and robbers. In this game, there are $s$ cops, and each stands on one vertex of the graph, and so does the robber. In each turn, first each cop moves to a neighbouring vertex, or does not move; and then the robber moves to a neighbouring vertex, or does not move. The cops win if at some stage, a cop is standing on the same vertex as the robber. Cops may share vertices and this game is played with full information. Given a graph $G$, how few cops suffice? The \emph{cop number} $c(G)$ is the smallest number of cops which can capture the robber $G$. The game played with one cop was initially defined by Nowakowski and Winkler \cite{nowakowski_vertex--vertex_1983} and Quilliot \cite{quilliot_these_1978}. The version with multiple cops was introduced Aigner and Fromme \cite{aigner_game_1984}; in particular they proved that the cop number of any connected planar graph is at most three.

Inspired by Andreae's result \cite{andreae_pursuit_1986} that the cop number of connected graphs forbidding $H$ as a minor is bounded for every graph $H$, Joret et al. \cite{joret_cops_2010} proved that the cop number of connected $H$-free graphs is bounded if and only if $H$ is a disjoint union of paths. In particular, they showed that the cop number of connected $P_t$-free graphs is at most $t-2$ for $t\geq 3$. Sivaraman \cite{sivaraman_application_2019} conjectured that two cops can win on connected $P_5$-free graphs and more generally that cop number of connected $P_t$-free graphs is at most $t-3$ for $t\geq 5$.

Other questions on the cop number and forbidden induced subgraphs have also been considered, for instance relating the cop number and the independence number (in other words, $tK_1$-free graphs) \cite{petr_note_2023} and excluding multiple induced subgraphs \cite{javadi_game_2022,masjoody_cops_2020,sivaraman_cop_2019-1}.

However, the question of the cop number of $P_5$-free graphs has received the most attention in this field; various weakenings of the conjecture about $P_5$-free graphs have been studied. Sivaraman and Testa \cite{sivaraman_cop_2019} conjectured the weaker statement that the cop number of connected $2K_2$-free graphs is at most two ($2K_2$, also written $2P_2$, is graph obtained by the disjoint union of two edges; it can also be seen as the complement of a four-vertex cycle). This was proved by the fourth author of this paper \cite{turcotte_cops_2022}. Liu \cite{liu_cop_2019} proved various partial results for these problems and Gupta, Mishra and Pradhan \cite{gupta_cops_2022} have proved the conjecture holds for multiple subclasses of $P_5$-free graphs. Masjoody \cite{masjoody_confining_2020} has conjectured the weaker statement that even if two cops perhaps cannot capture the robber on $P_5$-free graphs, they can confine it to a vertex.

In this paper, we prove Sivaraman's conjecture on the cop number of $P_5$-free graphs.

\begin{thm}\label{thm:cops}
	If $G$ is a connected $P_5$-free graph, then $c(G)\leq 2$.	
\end{thm}

The general strategy we employ is similar to the one used by the fourth author of this paper in \cite{turcotte_cops_2022} to bound the cop number of $2K_2$-free graphs. First, show that any graph in the class must contain a possible winning position for two cops, that is vertices $a,b\neq c$ such that $N[c]\subseteq N[a]\cup N[b]$. Then, consider a minimal graph in the class for which two cops cannot win, and use the minimality to force the robber to move to $c$, after which try to eventually move the cops to $a,b$ and show the robber cannot escape.

The first part of this strategy is accomplished by the following result, which we prove in \cref{sec:domineering}.

\begin{thm}\label{thm:domineering}
	If $G$ is connected and $P_5$-free, with $\alpha(G)\ge 3$, then there is a three-vertex induced path of $G$ with vertices $a, b, c$ in order, such that every neighbour of $c$ is also adjacent to one of $a,b$.
\end{thm}

Here, $\alpha(G)$ is the independence number of $G$, that is the cardinality of the largest stable subset of $V(G)$. Let $a\DD b \DD c$ be the vertices in order of a three-vertex induced path of $G$. We say that $a\DD b \DD c$ is {\em domineering}, or a {\em domineering 3-path},  if every neighbour of $c$ is also adjacent to one of $a,b$. Thus \ref{thm:domineering} says that every connected $P_5$-free graph with $\alpha(G)\ge 3$ has a domineering 3-path.

The condition $\alpha(G)\ge 3$ in \ref{thm:domineering} is needed. It is easy to see that every graph $G$ with $\alpha(G)=2$ is $P_5$-free, and it has no domineering 3-path if and only if its complement $H$ has diameter at most two, which gives plenty of counterexamples to \ref{thm:domineering} with $\alpha(G)\ge 3$ omitted.

The relative of \ref{thm:domineering} proved by the fourth author for $2K_2$-free graphs \cite{turcotte_cops_2022} is the following.
\begin{thm}\label{turcotte}
	If $G$ is connected and $2K_2$-free, with $|G|\ge 3$ and $G$ not a cycle of length five, then there exist distinct vertices $a,b,c$ such that $ab$ and $bc$ are edges (possibly $ac$ is also an edge) and every neighbour of $c$ is adjacent to one of $a,b$.
\end{thm}

This differs from \cref{thm:domineering} in three ways, two weakenings and a strengthening. First, it of course assumes that $G$ is $2K_2$-free, instead of $P_5$-free. Second, $ac$ might be an edge. But third, it does not need the assumption $\alpha(G)\ge 3$. What if we try to modify \ref{turcotte}, asking for a domineering 3-path in a $2K_2$-free graph $G$? Then again, it is false, but there are not so many counterexamples; every counterexample $G$ satisfies $\alpha(G)\le 2$, and it is easy to see that the counterexamples are the complements of Moore graphs of girth five. These are graphs of diameter two, with girth five, in which every vertex has the same degree $d$; such a graph exists only when $d=2,3,7$ and possibly $57$.

In order to prove \cref{thm:domineering}, we will need the following definition. Let us say a graph $G$ is {\em bijoined} if
\begin{itemize}[noitemsep]
	\item for every two nonadjacent vertices $u,v$ of $G$, there are exactly two vertices adjacent to both $u,v$, and they are adjacent to each other, and
	\item $G$ has no clique of cardinality four.
\end{itemize}
It is said to be {\em nontrivial} if $|G|>1$. If a nontrivial bijoined graph has a complement graph that is connected, then that complement would be a counterexample to \ref{thm:domineering}, so we care about bijoined graphs. Indeed, we will show in \cref{sec:domineering} that any counterexample to \ref{thm:domineering} has a connected induced subgraph whose complement is a nontrivial bijoined graph; so the whole question boils down to showing that there is no such graph. That is proved in \cref{sec:bijoined}.

The second part of the strategy, using \cref{thm:domineering} to capture the robber with two cops, is accomplished in \cref{sec:cops}. One important difference between the proofs for the $2K_2$-free case and the $P_5$-free case is that once the robber is on $c$, in the former  we can ensure the robber never leaves $c$, which is not possible in the latter.

Let us complete this section with some notation. Suppose $G$ is a graph, which we always consider to be simple and finite. For $v\in V(G)$, we write $N(v)$ for the \emph{neighbourhood} of $v$ (the set of vertices adjacent to $v$), $N[v]=N(v)\cup \{v\}$ for its \emph{closed neighbourhood} of $v$, and $M(v)=V(G)\setminus N[v]$ for the set of vertices distinct from and not adjacent to $v$. A vertex of a graph $G$ is {\em universal} if it is adjacent to every other vertex.

If $X,Y\subseteq V(G)$ are disjoint, we say $X$ is {\em complete} to $Y$ if every vertex in $X$ is adjacent to every vertex in $Y$, and $X$ is {\em anticomplete} to $Y$ if there are no edges between $X,Y$. If $X=\{x\}$, we say $x$ is {\em complete} to $Y$ if $\{x\}$ is complete to $Y$, and so on.

If $X\subseteq V(G)$, we write $G[X]$ for the subgraph of $G$ induced on $X$ and $G\setminus X$ for $G[V(G)\setminus X]$. If $X=\{v\}$, we write $G\setminus v$ for $G\setminus \{v\}$.

We will often write $x_1\DD x_2\DD \dots\DD x_k$ to represent a path with vertices $x_1, x_2, \dots, x_k$ in order.


\section{Bijoined graphs}\label{sec:bijoined}

We first note that bijoined graphs exist; for instance, if $H$ is a graph of girth at least five in which every two nonadjacent vertices have exactly one common neighbour,
and we add a universal vertex to $H$, we obtain a bijoined graph. 
We will show that no graphs are bijoined other than these, and in particular, no nontrivial bijoined graph has a connected 
complement graph.

We need the following well-known lemma:
\begin{thm}\label{moore}
Let $H$ be a graph with girth at least five,
such that for every two nonadjacent vertices $u,v$ there is exactly one vertex adjacent to both $u,v$. If the complement of $H$
is connected, then every two
vertices of $H$ have the same degree.
\end{thm}
\Proof
Since the complement of $H$ is connected, it suffices to show that every two nonadjacent vertices of $H$ have the same degree.
Thus, let $u,w$ be nonadjacent. 
Let $N=\{v_1\LL v_k\}$ be the set of neighbours of $u$, and for $1\le i\le k$ let $N_i$ be the
set of neighbours of $v_i$ different from $u$. Thus, the sets $\{u\},N, N_1\LL N_k$ are pairwise disjoint and have 
union $V(H)$. Let $w\in N_k$ say. For all $j\in \{1\LL k-1\}$, $w$ has a neighbour in $N_j$ (because it has 
distance two from $v_j$),
and has exactly one such neighbour (since $H$ has girth at least five); and has exactly one neighbour in $N\cup \{u\}$;
and so $w$ has degree exactly $k$. This proves \ref{moore}.~\bbox

We also need some results about strongly regular graphs. A graph is {\em strongly regular} with parameters $(n,k,a,c)$  
if it has $n$
vertices, every vertex has degree $k$, every two adjacent vertices have exactly $a$ common neighbours, and every two nonadjacent
vertices have exactly $c$ common neighbours.
Thus, Moore graphs are the strongly regular graphs that have parameters $(n,k,0,1)$ for some $n,k$. 
First, we need a result about Moore graphs 
mentioned earlier, due to Hoffman and Singleton~\cite{hoffman_moore_1960}:
\begin{thm}\label{mooregraphs}
A strongly regular graph with parameters $(n,k,0,1)$ 
exists only when $n=k^2+1$ and $k = 2,3,7$ and possibly 57.
\end{thm}
Second, we need the following (see Lemmas 10.3.2 and 10.3.3 of Godsil and Royle~\cite{godsil_algebraic_2001}):
\begin{thm}\label{godsilthm}
If a strongly regular
graph exists 
with parameters $(n,k,a,c)$, then 
either $2k=(n-1)(c-a)$ or $(a-c)^2+4(k-c)$ is a perfect square.
\end{thm}

Now we prove:
\begin{thm}\label{bijoin}
	If $G$ is a bijoined and non-null graph, then $G$ has a universal vertex.
\end{thm}
\Proof 
If $u,v$ are distinct, let $R(u,v)$ be the set of all vertices adjacent to both $u,v$; thus, if $u,v$ are nonadjacent
then $|R(u,v)|=2$
and its two members are adjacent. For convenience, we say ``$R(u,v)$ is an edge''.
We observe that no induced cycle of $G$ has length four (because for two opposite vertices $u,v$ of such a cycle,
$R(u,v)$ is not an edge). Thus $G$ is $C_4$-free, where $C_4$ is the cycle of length four. We suppose that $G$ has no 
universal vertex, for a contradiction. A {\em 4-clique} means a clique of cardinality four.

\begin{claim}\label{claim:bijoined:complementconnected}
	For each $v\in V(G)$, the complement of $G[N(v)]$ is connected.
\end{claim}

Suppose not. Then there is a partition $(X,Y)$ of $N(v)$ with $X,Y\ne \emptyset$, such that $X$ is complete to $Y$.
Since $G$ is $C_4$-free, one of $X,Y$ is a clique, say $X$; let $x\in X$. Then $x$ is adjacent to all other vertices in $N(v)$,
so we may assume that $X=\{x\}$ and $Y=N(v)\setminus \{x\}$. Since $G$ has no 4-clique, $Y$ is a stable set.
Since $x$ is not universal, there exists $z\in V(G)$ nonadjacent to $x$. Consequently $z\notin N[v]$. But then $R(v,z)$
is a subset of $Y$ and so not an edge, a contradiction. This proves \cref{claim:bijoined:complementconnected}.

\begin{claim}\label{claim:bijoined:samedegree}
	For each $v\in V(G)$, all vertices of $G[N(v)]$ have the same degree in $G[N(v)]$.
\end{claim}

$G[N(v)]$ has no cycle of length three, since $G$ is $K_4$-free; and $G[N(v)]$ is $C_4$-free since $G$ is $C_4$-free. Thus,
$G[N(v)]$ has girth at least five.
If $u,w\in N(v)$ are nonadjacent, then $R(u,w)$ consists of $v$ and exactly one vertex of $N(v)$; and so
$u,w$ have exactly one common neighbour in $G[N(v)]$. From \cref{claim:bijoined:complementconnected} and \ref{moore}, this proves \cref{claim:bijoined:samedegree}.

\bigskip

By \cref{claim:bijoined:samedegree} and since $G$ is bijoined, for every vertex $v$ there exists $k_v$ such that
$G[N(v)]$ is a Moore graph with parameters $((k_v)^2+1,k_v, 0, 1)$. In
particular, $|R(u,v)|=k_v$ for every neighbour $u$ of $v$.
It follows that $k_u=k_v$ for every pair of neighbours $u$ and $v$, and since $G$ is
connected, there exists $k$ such that $k_v=k$ for every $v$. So $G$ is a strongly
regular graph with parameters $(n,k^2+1,k,2)$ for some $n$ (which is determined
by $k$ but does not matter), where  $k \in \{2,3,7,57\}$ by \ref{mooregraphs}.
Let us apply \ref{godsilthm}, and deduce that one of the following holds:
\begin{itemize}
\item $2(k^2+1)=(n-1)(2-k)$; but this is impossible since $k\ge 2$.
\item $(k-2)^2+4((k^2+1)-2)= k(5k-4)$ is a perfect square; but this is not the case when $k=2,3,7,57$.
\end{itemize}
This contradiction proves \ref{bijoin}.~\bbox


\section{Finding a domineering 3-path}\label{sec:domineering}

Let us prove \ref{thm:domineering}, which we restate as follows.
\begin{thm}\label{thm:domineering2}
If $G$ is a connected $P_5$-free graph with $\alpha(G)\ge 3$, then there exists a domineering 3-path in $G$.
\end{thm}
\Proof
We assume that $G$ is a counterexample to the theorem with $G$ minimal. Thus, $G$ is connected and $P_5$-free, 
with $\alpha(G)\ge 3$, and there is no domineering 3-path in $G$, and no proper induced subgraph has these properties.
We will prove that the complement of $G$ is bijoined, which we will show is impossible. We begin with:

\begin{claim}\label{claim:domineering:nocorner}
	No two adjacent vertices $u,v$ satisfy $N(u)\subseteq N[v]$.
\end{claim}

Suppose that there are two such vertices $u,v$. If there is a vertex $w$ adjacent to $v$ and not to $u$, then $w\DD v\DD u$
is domineering, a contradiction; so $N[u]= N[v]$. Let $G'$ be obtained by deleting $u$. Then $G'$ is connected, $P_5$-free,
and satisfies $\alpha(G')\ge 3$, and so from the minimality of $G$, there is a domineering 3-path $a\DD b\DD c$ of $G'$.
This is not domineering in $G$, and so $u$ is adjacent to $c$ and nonadjacent to $a,b$. But then $v\ne a,b,c$, and so $v$
is adjacent to $c$ and not to $a,b$, contradicting that $a\DD b\DD c$ is domineering in $G'$. This proves \cref{claim:domineering:nocorner}.

\begin{claim}\label{claim:domineering:noncomplete}
	For each $v\in V(G)$ and  every component $C$ of $G[M(v)]$, no vertex $u\in N(v)$ is complete to $V(C)$.
\end{claim}

Because if $u$ is such a vertex, choose $w\in V(C)$; then $v\DD u\DD w$ is domineering, a contradiction. This proves \cref{claim:domineering:noncomplete}.

\begin{claim}\label{claim:domineering:connected}
	For each $v\in V(G)$, $G[M(v)]$ is non-null and connected. Moreover, every vertex in $N(v)$ has a neighbour in $M(v)$.
\end{claim}

By \cref{claim:domineering:nocorner}, $M(v)\ne \emptyset$. Suppose
that $C_1,C_2$ are distinct components of $G[M(v)]$. Since $G$ is connected, for $i=1,2$ there exists $u_i\in N(v)$
with a neighbour in $V(C_i)$. By \cref{claim:domineering:noncomplete}, for $i = 1,2$, $u_i$ has a neighbour and a non-neighbour in $V(C_i)$, and since $C_i$
is connected, there is an edge $a_ib_i$ of $C_i$ such that $u_i$ is adjacent to $a_i$ and not to $b_i$. If $u_1$ has a neighbour
in $V(C_2)$, we may assume that $u_1=u_2$, but then $b_1\DD a_1\DD u_1\DD a_2\DD b_2$ is a copy of $P_5$,
a contradiction. Thus $u_1$ has no neighbour in $V(C_2)$, and similarly $u_2$ has no neighbour in $V(C_1)$. If $u_1,u_2$
are nonadjacent then $b_1\DD a_1\DD u_1\DD v\DD u_2$ is a copy of $P_5$, and if $u_1,u_2$ are adjacent then
$b_1\DD a_1\DD u_1\DD u_2\DD a_2$ is a copy of $P_5$, in either case a contradiction. This proves the first assertion.
For the second, let $u\in N(v)$; then $u$ has a neighbour in $M(v)$ by \cref{claim:domineering:nocorner}. This proves \cref{claim:domineering:connected}.

\begin{claim}\label{claim:domineering:nonsubdomineering}
	For each $v\in V(G)$, $G[M(v)]$ has no domineering 3-path. Consequently $\alpha(G)=3$.
\end{claim}

Suppose that $a\DD b\DD c$ is a domineering 3-path of $G[M(v)]$. 
We claim that $a\DD b\DD c$ is also domineering in $G$. To show this, it suffices to show that
every neighbour $u$ of $c$ not in $M(v)$ is adjacent to one of $a,b$. But if not, then $a\DD b\DD c\DD u\DD v$
is a copy of $P_5$, a contradiction. This proves the first assertion of \cref{claim:domineering:nonsubdomineering}. For the second, suppose that
$\alpha(G)\ge 4$, and choose $v\in V(G)$ that belongs to a stable set of size four. Then $G[M(v)]$ has a stable
set of size three, and it is connected by \cref{claim:domineering:connected}, and has no domineering 3-path as we just showed, contrary to the minimality
of $|G|$. This proves \cref{claim:domineering:nonsubdomineering}.

\begin{claim}\label{claim:domineering:completeoranticomplete}
	For each edge $uv$, if $w\in N(u)\setminus N[v]$ and $C$ is a component of $G\setminus (N(u)\cup N(v))$, then $w$ is complete or anticomplete to $V(C)$.
\end{claim}

Suppose not; then there is an edge $ab$ of $C$ such that $w$ is adjacent to $a$ and not to $b$. But then
 $v\DD u\DD w\DD a\DD b$ is a copy of $P_5$, a contradiction. This proves \cref{claim:domineering:completeoranticomplete}.
 
\begin{claim}\label{claim:domineering:twocompletecomponents}
	If $v$ belongs to a stable set of size three, then for each edge $uv$, $G\setminus (N(u)\cup N(v))$ has exactly two components, both complete graphs.
\end{claim}

Certainly it has at most two components, since $\alpha(G)=3$, and for the same reason, if $G\setminus (N(u)\cup N(v))$
has two components then they are both complete graphs. Thus we just need to show that $G\setminus (N(u)\cup N(v))$ has at least two components.

Let $C= G\setminus (N(u)\cup N(v))$. Suppose that $C$ has at most one component.
By \cref{claim:domineering:nocorner}, $N(u)\setminus N[v]$ is nonempty. Let $w\in  N(u)\setminus N[v]$; then since $v\DD u\DD w$ is not domineering,
it follows that $w$ has a neighbour in $V(C)$, and hence is complete to $V(C)$ by \cref{claim:domineering:completeoranticomplete}. So $C$ is non-null, and 
$N(u)\setminus N[v]$ is complete to $V(C)$.
Now $M(v)= V(C)\cup (N(u)\setminus N[v])$, and $ V(C),N(u)\setminus N[v]$ are both nonempty.
If there is an induced path $a\DD b\DD c$ with $a,c\in N(u)\setminus N[v]$ and $b\in V(C)$, it follows that
$a\DD b\DD c$ is domineering in $G[M(v)]$, contrary to \cref{claim:domineering:nonsubdomineering}; and similarly there is no 
induced path $a\DD b\DD c$ with $b\in N(u)\setminus N[v]$ and $a,c\in V(C)$. Thus, $V(C)\cup (N(u)\setminus N[v])$ 
is a clique, contradicting that $v$ belongs to a stable set of size three. That proves \cref{claim:domineering:twocompletecomponents}.

\begin{claim}\label{claim:domineering:everyvertexstable3}
	Every vertex belongs to a stable set of size three; and so for every edge $uv$, $G\setminus (N(u)\cup N(v))$ has exactly two components, both complete graphs.
\end{claim}

Let $X$ be the union of all stable sets of size three. If $X\ne V(G)$, then since $G$ is connected, there is an edge
$uv$ with $u\notin X$ and $v\in X$. But then by \cref{claim:domineering:twocompletecomponents}, $G\setminus (N(u)\cup N(v))$ has  
exactly two components, both complete graphs, and consequently $u$ belongs to a stable set of size three, a contradiction.
This proves \cref{claim:domineering:everyvertexstable3}.

\begin{claim}\label{claim:domineering:nonadj}
	For every edge $uv$, $G\setminus (N(u)\cup N(v))$ consists of two nonadjacent vertices.
\end{claim}

By \cref{claim:domineering:everyvertexstable3}, $G\setminus (N(u)\cup N(v))$  has
exactly two components $C_1,C_2$, both complete graphs. For $i = 1,2$, let $X_i$
be the set of vertices in $N(u)\cup N(v)$ that have a neighbour in $V(C_i)$.
Suppose that $c_1,c_1'\in V(C_1)$ are distinct. From \cref{claim:domineering:everyvertexstable3} applied to the edge $c_1c_1'$, it follows that
$X_2\subseteq X_1$ (since otherwise the set of vertices nonadjacent to both $c_1,c_1'$ induces a connected subgraph).
Consequently, $N(u)\setminus N(v)$ is complete to $V(C_1)$.
Also, again by \cref{claim:domineering:everyvertexstable3} applied to the same edge, $N(u)\setminus N[v]\subseteq X_1$ and $N(v)\setminus N[u]\subseteq X_1$
(since for each $w\in  N(u)\setminus N[v]$, if $w\notin X_1$ then $\{v,u,w\}$ induces a three-vertex path, contrary
to \cref{claim:domineering:everyvertexstable3}).

Suppose that some $w\in N(u)\setminus N[v]$ belongs to $X_2$. Then $w$ is complete to $V(C_2)$ by \cref{claim:domineering:completeoranticomplete}, and the set of 
vertices nonadjacent to both $w,c_1$ is a subset of $N[v]$ including $v$ (because $w$ is complete to $V(C_1\cup C_2)$, 
and $c_1$ is complete to $N(u)\setminus N(v)$); and so this subset induces a connected subgraph, contrary to \cref{claim:domineering:everyvertexstable3}.
Thus $X_2\subseteq N(u)\cap N(v)$. If $c_2,c_2'\in V(C_2)$ are distinct, then the set of vertices nonadjacent to both $c_2,c_2'$
includes $u,v, w$ (where $w\in N(u)\setminus N[v]$), and these three vertices induce a path, contrary to \cref{claim:domineering:everyvertexstable3}. So $|C_2|=1$,
$C_2=\{c_2\}$ say. Choose $d\in N(u)\cap N(v)$ adjacent to $c_2$; then $u\DD d\DD c_2$ is domineering, a contradiction. This proves \cref{claim:domineering:nonadj}.

\bigskip

From \cref{claim:domineering:nonadj} and since $\alpha(G)=3$, it follows that the complement of $G$ is bijoined; and so from \ref{bijoin}, $G$ has a vertex
of degree zero, contradicting that it is connected. This proves \ref{thm:domineering2}.~\bbox


\section{Bounding the cop number}\label{sec:cops}

In order to prove our main result, we need the following definition and two lemmas regarding it.

We say that a subgraph $H$ of a graph $G$ is \emph{$P_3$-connected} if $H$ is connected, and for every pair of edges $e,f \in E(H)$ there exists a sequence of edges $e=e_0,e_1,\ldots,e_k=f$ such that $e_i$ and $e_{i+1}$ are two edges of an induced $P_3$ in $G$ for every $0 \leq i \leq k-1$. Note that the property of being $P_3$-connected is not an intrinsic property of $H$, but depends on $G$: even though $H$ is not required to be an induced subgraph of $G$, the pairs of edges in the definition must form induced paths in $G$, not just in $H$.

 \begin{thm}\label{thm:cops:anticompletepropagation}
 	If $H$ is a $P_3$-connected subgraph of a $P_5$-free graph $G$ and $u,v \in V(G)\setminus V(H)$ are such that $uv \in E(G)$, $u$ is anticomplete to $V(H)$ and the endpoints of some edge $e \in E(H)$ are non-neighbours of $v$, then $v$ is anticomplete to $V(H)$.
 \end{thm}
 
\Proof Suppose for a contradiction that $v$ has a neighbour in $V(H)$. As $H$ contains at least two vertices (as it contains $e$) and is connected, there exists $f \in E(H)$ such that its endpoint is a neighbour of $v$. Let $e=e_0,e_1,\ldots,e_k=f$ be the sequence of edges from the definition of $P_3$-connectedness (necessarily, $k\geq 1$). As $v$ is not adjacent to the endpoints of $e_0$ but is adjacent to at least one endpoint of $e_k$, there  exist $e_i,e_{i+1}$ such that $v$ is not adjacent to the endpoints of $e_i$ (say, $x,y$) but is adjacent to the other end of $e_{i+1}$ (say, $z$). By the choice of the sequence, $x$ is not adjacent to $z$. Then, $u\DD v\DD z\DD y\DD x$ is an induced $P_5$ in $G$, which is a contradiction. This completes the proof of \cref{thm:cops:anticompletepropagation}. ~\bbox

\begin{thm}\label{thm:cops:P3connectedexpansion} If $H$ is a $P_3$-connected subgraph of a connected graph $G$ and $v \in V(G)$, then either 
 \begin{enumerate}[noitemsep,label=\normalfont(\ref*{thm:cops:P3connectedexpansion}.\arabic*)]
 	\item\label{i:1} there exists a $P_3$-connected subgraph $H'$ of $G$ such that $H \subseteq H'$ and $v \in V(H')$, or
 	\item\label{i:2} there exists a $P_3$-connected subgraph $H'$ of $G$ such that $v \in V(H')$ and some $u \in V(H')$  is complete to $V(H)$ in $H'$, or
 	\item\label{i:3}  $v$ is complete to $V(H)$.
	\end{enumerate}
\end{thm}	

\Proof Let $Q$ be a shortest path in $G$ with one end $v$ and another end in $V(H)$. Let $v=v_0\DD v_1\DD \ldots\DD v_{\ell}$ be the vertex set of $Q$, where $v_{\ell} \in V(H)$. If $\ell = 0$, that is $v \in V(H)$, then $H'=H$ satisfies \ref{i:1}, and so we assume $\ell \geq 1$. 

Suppose now that $v_{\ell-1}$ has a non-neighbour in $V(H)$. Since $H$ is connected, we may suppose that $Q,v_\ell$ are chosen such that for some $w \in V(H)$ we have $v_{\ell}w \in E(H)$ and $v_{\ell-1}w \not\in E(G)$.  We claim that in this case $H' = Q \cup H$ satisfies \ref{i:1}; let us verify that $H'$ is  $P_3$-connected. As $Q$ is a shortest path, it is necessarily induced, and so $Q$ is $P_3$-connected. We also know that $H$ is $P_3$-connected. As $v_{\ell-1}v_{\ell} \in E(Q)$ and $v_{\ell} w \in E(H)$ are two edges of an induced $P_3$ in $G$ it follows that $Q \cup H$ is $P_3$-connected, as claimed.

It remains to consider the case when $v_{\ell-1}$ is complete to $V(H)$. If $\ell=1$, then $v=v_{\ell-1}$ and \ref{i:3} holds, so we may assume $\ell \geq 2$. Let $u=v_{\ell-1}$ and let $H'$ be a subgraph of $G$ with $V(H') = V(Q) \cup V(H)$ and $E(H')= E(Q) \cup \{uw : w \in V(H) \}$. We claim that $H'$ satisfies \ref{i:2}. Indeed, $u$ is complete to $V(H)$ in $H'$ and $H'$ is  $P_3$-connected since $Q$ is $P_3$-connected and every edge $uw \in E(H') - E(Q)$ forms a induced $P_3$ in $G$ with the edge $v_{\ell-2}u=v_{\ell-2}v_{\ell-1} \in E(Q)$.  This finishes the proof of \cref{thm:cops:P3connectedexpansion}.~\bbox

We are now ready to prove \cref{thm:cops}.

\bigskip
\Proof
	Suppose for a contradiction that there exists a connected $P_5$-free graph $G$ on which the robber has a winning strategy to evade two cops, and choose such $G$ with $|V(G)|$ minimum. In a series of claims we will obtain properties of this graph which will yield the desired contradiction.

\begin{claim}\label{claim:cops:indg2}
	$\alpha(G)>2$.
\end{claim}

To show \cref{claim:cops:indg2}, it suffices to see that if $\alpha(G)\leq 2$, there exists a dominating set of size at most two, on which two cops may start the game and win at the following turn, which is a contradiction.

\begin{claim}\label{claim:cops:existsdomineering}
	There exists a domineering path $a\DD b\DD r$ in $G$.
\end{claim}

This follows directly from \cref{thm:domineering} and \cref{claim:cops:indg2}. In the rest of the proof, we will always refer to a fixed domineering path $a\DD b\DD r$.

\begin{claim}\label{claim:cops:noncomplete}
	For every $v\in V(G)$, no vertex of $N(v)$ is complete to a component of $G[M(v)]$.
\end{claim}

Suppose otherwise that there exists $u\in N(v)$ and a component $C$ of $G[M(v)]$ such that $u$ is adjacent to every vertex of $C$. First note that $G\setminus C$ is connected as all vertices in $V(G)\setminus C$ adjacent to $C$ are in $N(v)$. By minimality of $G$, there exists a winning strategy for two cops on $G\setminus C$. We will use this strategy to define a winning strategy for two cops on $G$, which will be a contradiction. When playing on $G$, we say the robber's \emph{shadow} is on $x$ if the robber is on $x\in V(G)\setminus C$, and is on $v$ if the robber is on a vertex of $C$. In particular, the robber's shadow is always in $G\setminus C$. We show that any move of the robber yields a valid move for the robber's shadow in the sense that at every turn of the game the shadow either stays on its current vertex or moves to an adjacent vertex. At a given turn, suppose the robber moves from $x_1$ to $x_2$ (in particular, $x_1x_2\in E(G)$).
\begin{itemize}[noitemsep]
	\item If $x_1,x_2\in C$, then the robber's shadow stays on $v$, which is a valid move.
	\item If $x_1,x_2\notin C$, then the robber's shadow also moves from $x_1$ to $x_2$, which is a valid move
	\item If $x_1\in C$ and $x_2\notin C$, we note that necessarily $x_2\in N(v)$. Hence, the shadow moving from $v$ to $x_2$ is a valid move.
	\item If $x_1\notin C$ and $x_2\in C$, then $x_1\in N(v)$ and so the shadow moving from $x_1$ to $v$ is a valid move.
\end{itemize}
Consider the strategy for the cops on $G$ to follow the winning strategy on $G\setminus C$ to capture the robber's shadow. Once the cops have captured the robber's shadow, either they have captured the robber (if the robber and its shadow are on the same vertex) or the robber is in $C$ and its shadow is on $v$. In the latter case, having captured the shadow, at least one of the cops is on $v$. The other cop may then eventually move to $u$, and then capture the robber since $u$ is complete to $C$. Note that in the meantime, the robber cannot leave $C$ as it would be immediately captured by the cop on $v$. This proves \cref{claim:cops:noncomplete}.

Note that the proof of \cref{claim:cops:noncomplete} is a retract (special type of homomorphism) argument which is a standard tool in the study of the game of cops and robbers. A quite general version, which is close to the one presented here, was proved by Berarducci and Intriglia \cite{berarducci_cop_1993}.

\begin{claim}\label{claim:cops:nonneighbourdhoodconnected}
	For every $v\in V(G)$, $G[M(v)]$ is connected.
\end{claim}

Using \cref{claim:cops:noncomplete}, this is exactly the first part of \cref{claim:domineering:connected} in the proof of the existence of a domineering 3-path.

\begin{claim}\label{claim:cops:robberonor}
	There exists a strategy for two cops on $G$ to guarantee that, in order to avoid capture, the robber moves to $r$.
\end{claim}

By \cref{claim:cops:existsdomineering} we have that $N[r]\subseteq N[a]\cup N[b]$. Since $ab\in E(G)$, this implies that $G\setminus r$ is a connected graph. By minimality of $G$, there exists a winning strategy for two cops on $G\setminus r$. The strategy for cops in \cref{claim:cops:robberonor} is to play this strategy on $G$ as long as the robber has not entered $r$. If the robber never enters $r$ then it is eventually captured since this strategy is winning when restricted to $G\setminus r$. Hence, the robber eventually moves to $r$ (or chooses $r$ as its initial position). This proves \cref{claim:cops:robberonor}.

\vskip 10pt

In the rest of the proof, we will construct a strategy for two cops to attempt to capture the robber on $G$. However, since $c(G)>2$, we may assume that the robber has, and is playing, a strategy to avoid capture. The cops' strategy will begin by employing the strategy from \cref{claim:cops:robberonor}, and let $c_1,c_2 \in V(G)$ be the positions of the cops after the robber moves to $r$. Note that it is possible that $c_1=c_2$.

 \begin{claim}\label{claim:cops:instantwin}
 	$c_1,c_2 \in M(r)$.
 \end{claim}

 If \cref{claim:cops:instantwin} did not hold, one of the two cops could capture the robber at the next turn, which would be a contradiction.
 
 \begin{claim}\label{claim:cops:movingtodomineering}
 	For each $i\in \{1,2\}$, if $c_i \in N[a]$, then $c_{3-i} \not \in N(b)$.
 \end{claim}
 
 Suppose otherwise that for some $i\in \{1,2\}$, $c_i \in N[a]$ and $c_{3-i} \in N(b)$. Then the cops can move in one turn from $\{c_1,c_2\}$ to $\{a,b\}$. Being on $r$ at the start of the turn, the robber necessarily will be in $N[r]$ following its turn.  As $N[r]\subseteq N[a]\cup N[b]$ by \cref{claim:cops:existsdomineering}, the cops may then capture the robber at the following turn. This contradiction proves \cref{claim:cops:movingtodomineering}.
 
 \vskip 10pt
 
We say that a subgraph $H$ of $G$ is a \emph{snare} if 

 \begin{enumerate}[noitemsep,label=(H\arabic*)]
 	\item\label{item:HP3connected} $H$ is $P_3$-connected,
 	\item\label{item:Hnotadjr} $V(H) \subseteq M(r)$, 
 	\item\label{item:firstmove} there exists $d_1d_2 \in E(H)$ such that $d_1 \in N[c_1]$ and $d_2 \in N[c_2]$,
 	\item\label{item:ainH} $a \in V(H)$, and
 	\item\label{item:nonadjbinH} $V(H) \cap M(b) \neq \emptyset$.
 \end{enumerate}
 
 We finish the proof by showing that if a snare exists then the two cops can capture the robber, and finally that a snare exists.

\begin{claim}\label{claim:cops:snarethenwin}
	There is no snare.
\end{claim}

Suppose for a contradiction there exists a snare $H$. The cops are currently on $c_1,c_2$. The cops move to the endpoints of an edge $d_1d_2\in E(H)$ as in \ref{item:firstmove}. Since $H$ is connected, has at least one edge and $a\in V(H)$, there is an edge $f\in E(H)$ incident with $a$. By \ref{item:HP3connected}, there exists a sequence of edges $d_1d_2=e_0,e_1,e_2,\dots, e_k=f$ of $H$ such that $e_i$ and $e_{i+1}$ are two edges of an induced $P_3$ in $G$ for every $0 \leq i \leq k-1$. We may suppose without loss of generality that none of $e_0,\ldots e_{k-1}$ is incident with $a$. Over the next $k$ turns, the cops will follow this sequence of edges. In other words, in $i$ turns (for $0\leq i\leq k$) the cops will be on distinct endpoints of edge $e_i$. Finally, the cops will move to $\{a,b\}$. Note that by \ref{item:Hnotadjr}, $b \notin V(H)$. We remark that, since the robber is following a strategy which will avoid capture, it never enters the closed neighbourhood of a cop.

Let $r'$ be the first vertex visited by the robber in $M(b)$, and let $R$ be the set of all the previous positions of the robber. Note that $R\subseteq N[b]$. Such an $r'$ exists as once one of the cops is on $b$ (which always happens in the strategy described above), the robber must move to $M(b)$ as it would otherwise be captured at the next turn.

We claim that $R$ is anticomplete to $V(H)$. Suppose otherwise that at least one vertex of $R$ has a neighbour in $V(H)$. By \ref{item:Hnotadjr} $r$ has no neighbours in $V(H)$, and so there must exists consecutive positions of the robber $r_1,r_2 \in R$ such that $r_1$ has no neighbours in $V(H)$ (and is not itself in $V(H)$) but $r_2$ has at least one neighbour in $V(H)$. Since $r_2\in N[b]$, the cops are not on $\{a,b\}$ when the robber moves to $r_2$, and so they are positioned on distinct endpoints of an edge of $H$. These endpoints are non-neighbours of $r_2$, and so it follows from \cref{thm:cops:anticompletepropagation} that $r_2$ has no neighbours in $V(H)$. This is a contradiction, which implies our claim.

If cops are positioned on an edge of $H$ when the robber first moves to $r'$, then the same argument yields that $r'$ has no neighbour in $V(H)$, and because of \ref{item:ainH} we in particular have that $r'a\notin E(G)$. Otherwise, the cops are on $a,b$ when the robber moves to $r'$, and so $r'a \notin E(G)$. Hence, in all cases $r'a\notin E(G)$.

As $N[r]\subseteq N[a]\cup N[b]$, it follows that $r'r \not \in E(G)$. As $G[R]$ is connected, and $r'$ has a neighbour in $R$ there exists an induced path $r_1 \DD r_2 \DD r'$ such that $r_1,r_2 \in R$. Suppose that there exists an edge $xy \in E(H)$ such that $r'x \in E(G)$, but $r'y \notin E(G)$. Then $r_1\DD r_2\DD r'\DD x\DD y$ is an induced $P_5$. Hence, no such edge exists. As $H$ is connected and $r'a \not \in E(H)$ (in particular, $r'$ is not complete to $V(H)$) it follows that $r'$ is anticomplete to $V(H)$. Since $ab\in E(G)$ and $a\in V(H)$, it follows that $b$ is not anticomplete to $H$. However, by \ref{item:nonadjbinH} $b$ is not complete to $H$. Since $H$ is connected, there exists $xy \in E(H)$ such that $bx \in E(G)$ but $by \notin E(G)$. Then $r'\DD r_2\DD b\DD x\DD y$ is an induced $P_5$. This contradiction finishes the proof of \cref{claim:cops:snarethenwin}.

\vskip 10pt

It remains to show that there exists a snare in $G$, a contradiction.
	
\begin{claim}\label{claim:cops:snareexistence}
	There exists a snare.
\end{claim}	

Let $G'=G[M(r)]$. By \cref{claim:cops:nonneighbourdhoodconnected} $G'$ is connected. Let $P$ be an induced path in $G'$ with ends $c_1$ and $c_2$. Then $P$ is $P_3$-connected, and either $c_1=c_2$ or $P$ contains an edge satisfying \ref{item:firstmove}, as it has at most three edges.

Suppose first that $c_1 \in M(b)$. By \cref{thm:cops:P3connectedexpansion} applied to $G'$ with $H=P$ and $v=a$, either we find a subgraph $H'$ satisfying \ref{i:1} or \ref{i:2}, or $c_1,c_2 \in N[a]$. In the first and second cases, $H'$ is a snare;  property \ref{item:firstmove} is the only one which takes a little effort to verify. If $H'$ satisfies \ref{i:1} and $c_1\neq c_2$, then $P$ contains an edge satisfying \ref{item:firstmove} as noted above. If $H'$ satisfies \ref{i:1} and $c_1=c_2$, then $c_1x$ satisfies \ref{item:firstmove} for any neighbour $x$ of $c_1$ in $H'$; note that $x$ exists by connectivity of $H'$ unless $|V(H')|=1$, that is $a=c_1=c_2$, which contradicts that $c_1\in M(b)$. If $H'$ satisfies \ref{i:2} and $u \in V(H')$ is complete to $V(P)$ in $H'$ then $uc_1$ satisfies \ref{item:firstmove}. In the remaining case, $c_1,c_2 \in N[a]$ and $G'[\{c_1,a\}]$ is a snare. 

Thus we may assume that $c_1,c_2 \in N(b)$. By \cref{claim:cops:movingtodomineering} we have $c_1,c_2 \not \in N[a]$. By \cref{claim:cops:noncomplete} there exists $d \in V(G') \cap M(b)$.
Assume first that $\{c_1,c_2,d\}$ is not a clique of size three. By \cref{thm:cops:P3connectedexpansion} applied to $G'$ with $H=P$ and $v=d$, there exists a $P_3$-connected subgraph $H^*$ of $G'$ such that $c_1,c_2,d \in V(H^*)$, and $H^*$ contains a path with ends $c_1$ and $c_2$ with at most three edges. Indeed, if \ref{i:1} or  \ref{i:2} holds then $H^*=H'$ has the required properties. If  \ref{i:3} holds then $H^*=G'[\{c_1,c_2,d\}]$ is either an induced $P_2$ if $c_1=c_2$ or an induced $P_3$, as $d$ is adjacent to $c_1$ and $c_2$ and $\{c_1,c_2,d\}$ is not a clique of size three, and so $H^*$ is again as required.

We now apply \cref{thm:cops:P3connectedexpansion} to $G'$ with $H=H^*$ and $v=a$. As $c_1,c_2 \not \in N[a]$ and $c_1,c_2 \in V(H^*)$, $a$ is not complete to $V(H^*)$ and \ref{i:3} does not hold. Thus \ref{i:1} or  \ref{i:2} holds, that is there exists a $P_3$-connected subgraph $H'$ of $G'$ such that $c_1,c_2,d,a \in V(H')$ and $H'$ still contains  a path with ends $c_1$ and $c_2$ with at most three edges. It is routine to check that $H'$ is a snare.

It remains to consider the case when $C=\{c_1,c_2,d\}$ is a clique of size three.  We proceed similarly to the proof of \cref{thm:cops:P3connectedexpansion}. Let $Q$ be a shortest path from $a$ to $C$ with vertices $a=v_0\DD v_1\DD \ldots\DD v_{\ell}$, where $v_{\ell} \in C$. As $c_1,c_2 \not \in N[a]$, we have $a \not \in C$ and so $\ell \geq 1$.

Assume first that $\ell \geq 2$ and let $H$ be the subgraph of $G'$ defined by  $V(H)=V(Q) \cup C$, and 
$$E(H)= E(Q) \cup \{v_{\ell} u : u \in V(C), v_{\ell -1}u \not \in E(G')\} \cup \{v_{\ell-1} u : u \in V(C), v_{\ell -1}u  \in E(G')\}.$$
As $Q$ is $P_3$-connected, $H$ is connected, and every edge in $E(H) - E(Q)$ forms an induced $P_3$ in $G$ either with the edge $v_{\ell-2}v_{\ell-1}$ or $v_{\ell-1}v_{\ell}$, it follows that $H$ is $P_3$-connected. Note that $V(H) \subseteq V(G')=M(r)$ and $a, d \in V(H)$. Let us show with what choice of edge \ref{item:firstmove} holds. If at least one edge with both ends in $C$ is in $H$, then pick such an edge. Otherwise, it follows from the definition of $H$ that $v_{\ell-1}$ is necessarily complete to $C$, and so we can pick $c_1v_{\ell-1}$. Hence, $H$ is a snare.
 
It remains to consider the case $\ell=1$. In this case $a$ has a neighbour in $C$, and so $ad$ is the unique edge from $a$ to $C$, as $c_1,c_2 \not \in N[a]$. Then $H=G'[\{a,d,c_1\}]$ is an induced $P_3$ in $G'$ and the edge $dc_1$ satisfies  \ref{item:firstmove}. It follows that $H$ is a snare in this last case.

This completes the proof of \cref{claim:cops:snareexistence} and thus of the theorem.~\bbox

\section*{Acknowledgements}

This research was partially completed at the Second 2022 Barbados Graph Theory Workshop held at the Bellairs Research Institute in December 2022 and at the Combinatorics Workshop held at Mathematisches Forschungsinstitut Oberwolfach in January 2023.

\bibliography{refs}
\bibliographystyle{abbrvurl}

\end{document}